\documentclass[12pt]{article}
\usepackage{amsfonts}

 \def\baselinestretch{1} \topmargin -12pt \headsep 0pt \footskip 30pt
 \def\medskipamount{12pt} \def\smallskipamount{6pt}

\newcounter{bitcount}
\newcommand{\bit}[1]{\addtocounter{bitcount}{1}\pagebreak[3]
\subsection{#1}\nopagebreak\setcounter{equation}{0}}

\renewcommand{\theequation}{\thesubsection .\arabic{equation}}
\renewcommand{\thesubsection}{\arabic{bitcount}}

\catcode`\@=\active
\catcode`\@=11
\def\@eqnnum{\hbox to .01pt{}\rlap{\bf \hskip -\displaywidth(\theequation)}}
\catcode`\@=12

\begin{document}


\catcode`\@=\active
\catcode`\@=11
\newcommand{\nc}{\newcommand}


\nc{\bs}[1]{ \addvspace{\medskipamount} \pagebreak[3]
\refstepcounter{equation}
\noindent {\bf (\theequation) #1.} \begin{em} \nopagebreak}

\nc{\es}{\end{em} \par\addvspace{\medskipamount} } 

\nc{\br}[1]{ \addvspace{\medskipamount} \pagebreak[3]
\refstepcounter{equation} 
\noindent {\bf (\theequation) #1.} \nopagebreak}

\nc{\er}{\par \addvspace{\medskipamount} }


\nc{\C}{\mathbb C}
\nc{\Pj}{\mathbb P}
\nc{\Q}{\mathbb Q}
\nc{\R}{\mathbb R}
\nc{\Z}{\mathbb Z}


\nc{\oper}[1]{\mathop{\mathchoice{\mbox{\rm #1}}{\mbox{\rm #1}}
{\mbox{\rm \scriptsize #1}}{\mbox{\rm \tiny #1}}}\nolimits}
\nc{\ad}{\oper{ad}}
\nc{\diag}{\oper{diag}}
\nc{\imag}{\oper{Im}}
\nc{\tr}{\oper{tr}}

\nc{\operlim}[1]{\mathop{\mathchoice{\mbox{\rm #1}}{\mbox{\rm #1}}
{\mbox{\rm \scriptsize #1}}{\mbox{\rm \tiny #1}}}}

\nc{\Ui}{{{\rm U(1)}}}
\nc{\Uii}{{{\rm U(2)}}}
\nc{\Suii}{{{\rm SU(2)}}}
\nc{\suii}{{{\mathfrak s \mathfrak u{\rm(2)}}}}
\nc{\ui}{{{\mathfrak u{\rm(1)}}}}

\nc{\Ga}{\Gamma}
\nc{\ga}{\gamma}
\nc{\La}{\Lambda}
\nc{\la}{\lambda}


\nc{\drow}{\raisebox{.4ex}[3ex]{\hbox{$\left\downarrow\vbox to
    9.5pt{}\right.\nulldelimiterspace=0pt \mathsurround=0pt$}}}



\nc{\ratio}[2]{\mathchoice{ {\textstyle \frac{\scriptstyle
#1}{\scriptstyle #2}}} {{\textstyle\frac{\scriptstyle #1}{\scriptstyle
#2}}} {\frac{\scriptscriptstyle #1}{\scriptscriptstyle #2}}
{\frac{\scriptscriptstyle #1}{\scriptscriptstyle #2}}}

\nc{\half}{\mathchoice{{\textstyle \frac{\scriptstyle 1}{\scriptstyle
2}}} {{\textstyle\frac{\scriptstyle 1}{\scriptstyle 2}}}
{\frac{\scriptscriptstyle 1}{\scriptscriptstyle 2}}
{\frac{\scriptscriptstyle 1}{\scriptscriptstyle 2}}}


\nc{\V}{V}
\nc{\cl}{\Gamma}
\nc{\rs}{X}


\nc{\Left}[1]{\hbox{$\left#1\vbox to
    11.5pt{}\right.\nulldelimiterspace=0pt \mathsurround=0pt$}}
\nc{\Right}[1]{\hbox{$\left.\vbox to
    11.5pt{}\right#1\nulldelimiterspace=0pt \mathsurround=0pt$}}


\nc{\down}{\Big\downarrow}
\nc{\fp}{\mbox{     $\maltese$} \par \addvspace{\smallskipamount}}
\nc{\lrow}{\longrightarrow}
\nc{\pf}{\noindent {\em Proof}}
\nc{\st}{\, | \,}
\nc{\sans}{\backslash}

\catcode`\@=12



\noindent
{\LARGE \bf A perfect Morse function\\ 
on the moduli space of flat connections}
\medskip \\ 
{\bf Michael Thaddeus } \smallskip \\ 
Department of Mathematics, Columbia University \\
2990 Broadway, New York, N.Y. 10027
\renewcommand{\thefootnote}{}
\footnotetext{Partially supported by NSF grants DMS--9500964 and
  DMS--9808529.}

\bigskip \bigskip

\noindent
Let $\rs$ be a compact surface of genus $g$, and let $M_g$ be
the moduli space of flat $\Suii$ connections on $\rs$ having holonomy
$-I$ around a single puncture $p$.  Let $a_1, b_1, a_2, b_2, \dots,
a_g, b_g$ be the usual generators for $\pi_1(\rs \sans p)$, and define
a real-valued function on $M_g$ by assigning to a flat connection the
trace of its holonomy around $a_g$.  This paper will give an
elementary, direct proof that this function is a {\em perfect}
Bott-Morse function, that is, one whose Morse inequalities are
equalities.  This leads to a new derivation of the well-known
Harder-Narasimhan formula for the Betti numbers of the moduli space.

The Harder-Narasimhan formula has previously been calculated by
several methods.  Most of them involved the deep theorem of Narasimhan
and Seshadri \cite{ns} identifying flat connections with stable
holomorphic bundles.  For example, Harder and Narasimhan \cite{hn}
applied the Weil conjectures, counting algebraic vector bundles on
curves over finite fields to deduce the result.  Another approach, due
to Atiyah and Bott \cite{ab}, used two-dimensional Yang-Mills theory
to construct $M_g$ as the symplectic quotient of an
infinite-dimensional affine space by the gauge group.

More recently, the Bott-Morse function discussed here was
studied by Jeffrey and Weitsman.  It is straightforward to determine
the critical points of this function (cf.\ \ref{crit}), and hence to
compute one side of the Morse inequalities.  Jeffrey and Weitsman
noted that the expression thus obtained equals the Harder-Narasimhan
formula.  If the latter is regarded as known, this constitutes an {\em
a posteriori} proof that the function is perfect.

This paper will give an {\em a priori} proof that the function is
perfect, using only relatively elementary facts from
finite-dimensional symplectic geometry.  It therefore provides a
down-to-earth proof of the Harder-Narasimhan formula, avoiding high
technology such as the Weil conjectures or Yang-Mills theory.
It is much closer in spirit to the original computation of the Betti
numbers, due to Newstead \cite{h, n}.

A few years ago, Mark Hoyle pointed out to the author that the {\em a
posteriori} argument remains valid for the space of flat $\Suii$
connections whose holonomy around $p$ lies in any fixed conjugacy
class (except that of $I$ when there are reducibles).  In this case
the Betti numbers are easily determined from the Harder-Narasimhan
formula. The present {\em a priori} argument extends without change to
this space.  Even better, it applies more generally to the space of
flat $\Suii$ connections on $\rs$ with any finite number of punctures,
provided only that the conjugacy classes of the holonomies around the
punctures are fixed so that there are no reducibles.  Meanwhile Hoyle
has obtained formulas for the Betti numbers in this case too and shown
that the {\em a posteriori} argument extends \cite{hoy}.  The results
of this paper therefore provide an independent proof of his formulas.

A natural question for further investigation is how much can be said
for other Lie groups than $\Suii$.  One could replace the trace with
any real-valued smooth class function.  The proof of the main theorem
does not carry over in full, but it might imply some partial results.
Furthermore, the {\em a posteriori} question should be tractable.

The organization of the paper is simple.  Section \ref{bkg} reviews
what little background material is needed: the construction of the
moduli space as a space of representations of $\pi_1(X)$, and some
results of Goldman and Frankel.  The main theorem is stated, and the
Poincar\'e polynomial of the moduli space is derived.  Section
\ref{spec} describes the spectral sequence, canonically associated to
a Bott-Morse function $f: M \to \R$, which abuts to $H^*(M, \Q)$.
Section \ref{pf} is devoted to the proof of the main theorem.  The
basic strategy is to prove that all differentials vanish in the
spectral sequence by combining the perfection of moment maps with a
parity argument.  Section \ref{punct} explains how to generalize the
main theorem to the case of a surface with many punctures, recovering
the formulas of Hoyle.  Finally, \S\ref{action} explains why the
methods of \S\S\ref{spec} and \ref{pf} immediately show that the
integral cohomology is torsion-free, and yield information about the
$\Uii$ case as well.

Unless otherwise specified, all cohomology is with rational
coefficients. 

This work was announced in lectures given in Odense and Barcelona
\cite{thad}.  I wish to apologize for the delay in publication.

I am very grateful to Kenji Fukaya, Mark Hoyle, Michael Hutchings, Lisa
Jeffrey, Tom Mrowka, and Jonathan Weitsman for helpful discussions.
I especially thank Hoyle for providing me with a copy of
his paper \cite{hoy} in advance of publication.

\bit{The moduli space of flat \boldmath$\Suii$ connections}
\label{bkg}

The moduli space $M_g$ which interests us can be defined simply as
follows.  Let $\mu_g: \Suii^{2g} \to \Suii$ be given by $(A_1, B_1,
\dots, A_g, B_g) \mapsto - \prod_i [A_i, B_i]$.  Then $I$ is a regular
value of $\mu_g$.  (This was first shown by Igusa \cite{igu}; it is
also a special case of Proposition \ref{critval} below.)  Hence
$\mu_g^{-1}(I)$ is smooth; moreover $\Suii / \pm I$ acts freely on it
by conjugation, so $M_g = \mu_g^{-1}(I) \, / \, \Suii$ is a smooth
compact $6g-6$-dimensional manifold.

At any point $\rho \in M_g$, the tangent space is naturally isomorphic
to the first cohomology of the complex
$$\suii \lrow \suii^{2g} \lrow \suii, $$
where the former map is the derivative of conjugation at $(A_i, B_i)$,
and the latter is the derivative of $\mu_g$.  
Choose any representative for $\rho$ in $\mu_g^{-1}(I)$ and, by abuse
of notation, denote it again by $\rho$.  Then 
$\ad \rho$ defines a representation of $\pi_1(\rs)$ on $\suii$, and
the cohomology mentioned above is none other than the group cohomology
$H^1(\pi_1(X), \, \ad \rho)$. 

Since $\rs$ is an Eilenberg-Mac\phantom{.}Lane space, $H^2(\pi_1(\rs),
\R) = H^2(\rs, \R) = \R$.  Combining the cup product with the
symmetric form $\langle a, b \rangle = - \half \tr ab$ gives
a non-degenerate antisymmetric map
$$ H^1(\ad \rho) \otimes H^1(\ad \rho)
\stackrel{\cup}{\lrow} H^2(\ad\rho \otimes \ad \rho) 
\stackrel{\langle\,\, , \,\,\rangle \hspace{0.2em}}{\lrow} H^2(\R) = \R, $$
which determines a non-degenerate 2-form on $M_g$.

\bs{Theorem (Goldman)}
This 2-form is closed.
\es

Thus $M_g$ becomes a symplectic manifold.  Goldman's original
proof of this theorem \cite{gold1} used infinite-dimensional quotients
in the style of Atiyah and Bott.  Since then, purely algebraic,
finite-dimensional proofs have been provided by Karshon \cite{kar} and
Weinstein \cite{wein}; see also Guruprasad et al.\ \cite{guru}.  

Now let $f: M_g \to [-1,1]$ be given by $(A_i, B_i) \mapsto \half \tr
A_g$, which is well-defined since the trace is conjugation-invariant.
Then $\Ui$ acts on $f^{-1}(-1,1)$ as follows.  If $A_g \neq \pm I$,
then there is a unique homomorphism $\phi: \Ui \to \Suii$ such that
$A_g \in \phi(\{ \imag z > 0 \} )$.  Let $\la \cdot (A_i, B_i) = (A_1,
B_1, \dots, A_g, B_g \cdot \phi(\la))$.  

\bs{Proposition (Goldman)}
\label{gold}
This action preserves the symplectic form, and it has moment map
$-i\arccos f$.
\es

See Kirwan \cite{k} for the definition of a moment map, and Goldman
\cite{gold2} and Jeffrey-Weitsman \cite{jw2} for a proof.

If the $\Ui$-action were global on $M_g$, the following result
\cite{frank} would immediately imply that $\arccos f$ was perfect.

\bs{Theorem (Frankel)}\label{kir}%
The moment map of a $\Ui$-action on a compact
symplectic manifold (times $i$) is a
perfect Bott-Morse function.
\es

In the present case, however, the $\Ui$-action on $f^{-1}(-1,1)$ does
not extend over $f^{-1}(\pm 1)$, and $\arccos f$ is not even
differentiable there.  Nevertheless, the following theorem is true,
and will be proved in \S\ref{pf}.

\bs{Main Theorem}
The map $f$ is a perfect Bott-Morse function on $M_g$.
\es

To see what this means concretely, let us identify the critical
submanifolds and their indices, assuming for the moment that $f$ is a
Bott-Morse function.  

\br{The critical submanifolds of \boldmath$f$}
\label{crit}
First, let $S_1 = f^{-1}(-1)$.  As the absolute minimum of $f$, $S_1$
is of course a critical submanifold.  It is exactly the locus where
$A_g = -I$; hence $B_g$ may be arbitrary, and the product of the first
$g-1$ commutators must be $-I$, so $S_1 = (\mu_{g-1}^{-1}(I) \times
\Suii) / \Suii$.  The natural projection $\pi: S_1 \to M_{g-1}$ makes
$S_1$ into an $\Suii$-bundle over $M_{g-1}$.  This
is an adjoint, not a principal bundle, so it may have a section
without being trivial.  Indeed, $B_g = I$ determines such a section;
hence the Euler class vanishes and so by the Gysin sequence 
$H^*(S_1) = V \oplus \tau V$, where $V = \pi^* H^*(M_{g-1})$ and $\tau
\in H^3(S_1)$ is the Poincar\'e dual of the locus where $B_g = I$. 
As the absolute minimum, $S_1$ of course has index 0.

Exactly the same is true of $S_3 = f^{-1}(1)$, except that it is
the absolute maximum of $f$, the locus where $A_g = I$, and hence has
index equal to its codimension, which is 3.

Within $f^{-1}(-1, 1)$, on the other hand, the critical points of $f$
coincide with those of $\arccos f$, the moment
map for the $\Ui$-action. They are therefore exactly the fixed points
of that action, and hence are represented by $2g$-tuples $(A_i, B_i)
\in \Suii^{2g}$ that are conjugate to $(A_1, B_1, \dots, A_g, B_g
\cdot \phi(\la))$ for all $\la \in \Ui$.  It is straightforward to
check that these are all conjugate to $2g$-tuples such that 
$$A_g = \left(\begin{array}{cr}i&0\\0&-i\end{array}\right), \, \, \,
B_g = \left(\begin{array}{rc}0&1\\-1&0\end{array}\right),$$ 
and the remaining $A_i$ and $B_i$ are diagonal.
Hence the only other critical value is $\half \tr A_g$, which is
0, and the corresponding critical set $S_2$ is a $2g-2$-torus.
Because there is an involution $A_g \mapsto -A_g$ on $M_g$ changing
the sign of $f$, the index of $S_2$ must be half the rank of its
normal bundle, or $2g-2$. 

Incidentally, $S_1$ and $S_3$ are empty when $g=1$, and everything is
empty when $g=0$, but these special cases will not affect our arguments.
\er

The main theorem therefore implies the so-called Harder-Narasimhan
formula. 

\bs{Corollary} The Poincar\' e polynomial of $M_g$ is
$$P_t(M_g) = \frac{(1+t^3)^{2g} - t^{2g}(1+t)^{2g}}{(1-t^2)(1-t^4)}.$$
\es \vspace{-3 ex}

\pf.  Since $M_0$ is empty, $P_t(M_0) = 0$ as desired.  It follows
from the theorem and the discussion above that
\begin{eqnarray*}
P_t(M_g) & = & P_t(S_1) + t^3 P_t(S_3) + t^{2g-2} P_t(S_2) \\
& = & (1+t^3)^2 P_t(M_{g-1}) + t^{2g-2}(1+t)^{2g-2}.\end{eqnarray*}
This gives a recursion for $P_t(M_g)$, which is satisfied by the
formula.  \fp

\bit{The spectral sequence of a Bott-Morse function}
\label{spec}

Let $M$ be a compact manifold and $f: M \to \R$ a {\em Bott-Morse
function}, that is, a smooth function whose critical set is a
disjoint union of submanifolds on whose normal bundles the Hessian is
nondegenerate.  For $f$ to be perfect is clearly equivalent to the
vanishing of all differentials in the spectral sequence described
below. 

The set of critical values of $f$ is finite, so the critical set of
$f$ is a disjoint union of submanifolds $S_1, S_2, \dots, S_n \subset
M$ such that $f$ is constant on each $S_i$ and $f(S_i) < f(S_j)$ for
$i < j$.  Choose $x_j \in \R$ satisfying $x_0 < f(S_1) < x_1 < f(S_2)
< x_2 < \cdots < f(S_n) < x_n$ and let $U_j = f^{-1}{(x_0,x_j)}$, so
that $M$ is filtered by open sets $\emptyset = U_0 \subset U_1 \subset
\cdots \subset U_n = M$.  The complex of singular chains on $M$ is
then filtered by the support of the chain:
$$0 = C_*(U_0) \subset C_*(U_1) \subset \cdots \subset C_*(U_n) =
C_*(M).$$
Moreover, this filtration is preserved by the differential.  

Taking duals yields a cofiltration on the group of cochains, which is
also preserved by the differential.  This is exactly the raw material
needed to construct a spectral sequence, and the whole machine runs
smoothly.  Except that everything is dualized because of the
cofiltration, it is just as described by Bott \& Tu \cite{bt}.  As in
their equation (14.3), there is a short exact sequence
$$0 \lrow B \lrow A \stackrel{\pi}{\lrow} A \lrow 0$$
where $A = \bigoplus_j C^*(U_j)$ and $\pi: C^*(U_j) \to C^*(U_{j-1})$ is
the restriction of cochains.  This leads to an exact couple
$$\begin{array}{c}
\bigoplus H^*(U_j) \lrow \bigoplus H^*(U_j) \\
\nwarrow \phantom{XXXXX} \swarrow \\
\bigoplus H^*(U_j,U_{j-1})
\end{array}$$
whose derived couples abut to $H^*(M)$.  On the other hand, the Morse
lemma implies that up to homotopy, $U_j$ is a CW complex obtained
from $U_{j-1}$ by attaching, along its boundary, the disc bundle $E_j$
associated to the {\em negative normal bundle} consisting of the
negative definite subspaces of the Hessian of $f$ on $S_j$.  So by
excision, $H^*(U_j, U_{j-1}) \cong H^*(E_j, \partial E_j)$.

On the other hand, suppose that $S_j$ and $E_j$ are orientable for
each $j$, and choose orientations.  Then the cup product with the
Thom class induces the {\em Thom isomorphism} $H^*(E_j, \partial E_j)
\cong H^*(S_j)$.  The spectral sequence then can be regarded as a
sequence of differentials $d_1, d_2, \dots$ where $d_1: \bigoplus
H^*(S_j) \to \bigoplus H^*(S_j)$ and $d_{i+1}: H^*(d_i) \to H^*(d_i)$.

By the definition of the differentials of an exact couple, the first
differential $d_1$ is the direct sum of maps $H^*(S_j)
\to H^*(S_{j+1})$ induced by the upper part of the following diagram.  

{\def\arraycolsep{1pt}
\begin{equation}\label{diag}\begin{array}{crcccccccc}
& H^*(S_j) &&&&&&&& \\
&& \searrow &&&&&&& \\
&&& \multicolumn{3}{c}{H^*(U_j, U_{j-1})} &&&& \\
&&&& \drow &&&&& \\
& H^*(U_{j+1}) & \multicolumn{2}{c}{\lrow} & H^*(U_j) &
\multicolumn{2}{c}{\lrow} & H^*(U_{j+1}, U_j) && \\
&\multicolumn{1}{c}{\drow} &&&&&&& \raisebox{.4ex}{$\searrow$} & \\
\multicolumn{2}{c}{H^*(U_{j+2}, U_{j+1})} &&&&&&&& H^*(S_{j+1}) \\
&& \searrow &&&&&&& \\
&&&& \multicolumn{1}{l}{H^*(S_{j+2})} &&&&&
\end{array}\end{equation}}%
Here all the maps are induced by inclusion, except the diagonal
arrows, which are Thom isomorphisms.  

An element of $H^*(S_j)$ is thus in the kernel of the first
differential if and only if it maps to $0 \in H^*(U_{j+1}, U_j)$.  By
exactness of the middle row, it comes from some element of
$H^*(U_{j+1})$, and its image in $H^*(S_{j+2})$ is the value of the
second differential.

The higher differentials can now be described in the same manner, 
applying the argument of the previous paragraph inductively.

The contents of this section go back to Bott \cite{b}.

\bit{Proof of the main theorem}%
\label{pf}

Let $M_g$ and $f: M_g \to \R$ be as in \S\ref{bkg}.  To prove the main
theorem, we will show first that $f$ is a Bott-Morse function, then
that it is perfect.  The former task is accomplished with three
lemmas; the latter will occupy the remainder of the section.

\bs{Lemma}%
\label{lemma}%
For $g \geq 2$, the map $\mu_g^{-1}(I) \to \Suii \times \Suii$ given
by $(A_i, B_i) \mapsto (A_g, B_g)$ is a submersion at the locus where
$A_g = \pm I$. \es

\pf.  This means that the infinitesimal map $\ker D \mu_g \to \suii
\times \suii$ induced by projection of $\Suii^{2g}$ on the last two
factors is surjective when $A_g = \pm I$.  Direct computation shows
that in this case,
$$D \mu_g(a_1, \dots, b_g) = 
D\mu_{g-1}(a_1, \dots, b_{g-1}) + D\mu_1(a_g, b_g).$$
But as mentioned before, $I$ is a regular value of
$\mu_{g-1}$.  Hence $D \mu_{g-1}$ is surjective, so there exist
$(a_i, b_i) \in \ker D \mu_g$ having any desired values for $a_g$ and
$b_g$.  \fp

\bs{Lemma}%
\label{lm}%
There is an $\Suii$-equivariant diffeomorphism between
a neigh\-bor\-hood of $S_1 \subset \mu_g^{-1}(I)$ and a
neighborhood of 
$$\mu_{g-1}^{-1}(I) \times \{ -I \} \times \Suii
\subset \mu_{g-1}^{-1}(I) \times \Suii \times \Suii,$$ 
identifying the map of the previous lemma with the projection on the
last two factors.  Likewise for $S_3$ if $\{ I \}$ is substituted for
$\{ -I \}$. \es

\pf.  This follows immediately from the previous lemma using the
equivariant version of the tubular neighborhood theorem for
submanifolds \cite[I Thm.\ 2.1.1]{audin} and the inverse function
theorem. \fp 

\bs{Lemma}
The function $f$ is a Bott-Morse function.
\es

\pf.  By definition, we need to show that the Hessian on the normal
bundle to every critical submanifold is nondegenerate.  By Frankel's
result, Theorem \ref{kir}, $\arccos f$ is a Bott-Morse function on
$f^{-1}(-1,1)$, and hence so is $f$.  Therefore the Hessian of $S_2$
is nondegenerate.  In a neighborhood of $S_1$, $f$ is locally the
composite with the diffeomorphism from Lemma \ref{lm} of the trace map
on the first $\Suii$ factor.  But if $\Suii$ is identified with the
3-sphere of unit quaternions, then the trace map is simply a linear
projection, whose Hessian at $\pm I$ is certainly nondegenerate.
Hence the Hessian of $f$ is nondegenerate on $S_1$.  The case of $S_3$
is similar.  \fp
 
\br{The symplectic cut \boldmath$\hat{M}_g$} It remains to show that $f$
is perfect.  As a first step towards this goal, we construct a slight
modification of $M_g$ on which the $\Ui$-action becomes global and hence
the moment map is perfect.  Take the 2-sphere $S^2$ with the standard
$\Ui$-action given by rotation, and normalize the moment map so that
its image is $[-\ratio{2}{3}, \ratio{2}{3}]$.  Then let $\hat{M}_g$ be
the symplectic quotient of $f^{-1}(-1,1) \times S^2$ by the diagonal
$\Ui$-action.  This is a simple example of a {\em symplectic cut\/} of
$M_g$, as introduced by Lerman \cite{ler}.  As a set, it is just
$f^{-1}[-\ratio{2}{3}, \ratio{2}{3}]$ with the $\Ui$-orbits in
$f^{-1}(-\ratio{2}{3})$ and $f^{-1}(\ratio{2}{3})$ collapsed to
points.  However, the symplectic cut construction shows that it is a
symplectic manifold and that the residual $\Ui$-action is Hamiltonian.
Since this action is globally defined, its moment map $\hat{f}:
\hat{M}_g \to [-\ratio{2}{3}, \ratio{2}{3}]$ is perfect.  It again has
three critical submanifolds $\hat{S}_1, \hat{S}_2,$ and $\hat{S}_3$.  Of
these, $\hat{S}_2 = S_2$ obviously.  However, $\hat{S}_1$ and
$\hat{S}_3$ are $S^2 \times S^2$-bundles over $M_{g-1}$; this follows
from Lemma \ref{lm}. \er

Filter $M_g$ by the open subsets 
$U_0 = f^{-1}(-\ratio32, -\ratio32)$,
$U_1 = f^{-1}(-\ratio32, -\ratio12)$, 
$U_2 = f^{-1}(-\ratio32, \ratio12)$, and
$U_3 = f^{-1}(-\ratio32, \ratio32)$, so that
$$\emptyset = U_0 \subset U_1 \subset U_2 \subset U_3 = M_g.$$  
The machinery of \S\ref{spec} then
produces a spectral sequence abutting to $H^*(M_g)$.
We need to show that all the differentials are zero.  
Since there are only three nonzero terms in the cochain filtration, 
only the first two differentials can possibly be nonzero.  

As shown in \S\ref{spec}, the first differential is the direct sum of
the maps $H^*(S_1) \to H^*(S_2)$ and $H^*(S_2) \to H^*(S_3)$ induced
by the upper route in the diagram (\ref{diag}).  We will consider these
maps separately.

\br{The map \boldmath $H^*(S_1) \to H^*(S_2)$} 
\label{s1s2}
As seen in (\ref{crit}),
$H^*(S_1) = \V \oplus \tau \V$, where $\tau \in H^3$ is Poincar\'e
dual to the locus where $B_g = I$, and $\V = \pi^* H^*(M_{g-1})$.  The
map $H^*(S_1) \to H^*(S_2)$ accordingly splits into two maps $\V \to
H^*(S_2)$ and $\tau \V \to H^*(S_2)$; we will show that each of these
vanishes.

Let $U'_1 = f^{-1}(-\ratio23,-\ratio12)$.  Then $U'_1 \subset U_1$,
but there is also a natural inclusion $U'_1 \subset \hat{U}_1$, where
$\hat{U}_1 = \hat{f}^{-1}(-\ratio32, -\ratio12)$.  It
then follows from Lemma \ref{lm} that the hexagon in the
diagram below commutes.  Here the maps of the second column are the Thom
isomorphisms, which in this case are simply induced by retractions,
and those of the first column are induced by the fiber bundle
projections.
{\def\arraycolsep{1pt}
$$\begin{array}{ccccccccccc}
&& H^*(S_1) & \lrow & H^*(U_1) &&&&&& \\
& \nearrow &&&& \searrow &&&&& \\
H^*(M_{g-1}) &&&&&& H^*(U'_1)& \lrow & H^*(U'_2, U'_1) & \lrow &
H^*(S_2) \\ 
& \searrow &&&& \nearrow &&&&& \\
&& H^*(\hat{S}_1) & \lrow & H^*(\hat{U}_1) &&&&&& 
\end{array}$$}

On the other hand, if $U'_2 = f^{-1}(-\ratio23, \ratio12)$, then by
excision $H^*(U'_2, U'_1) = H^*(U_2, U_1)$.  The map $\V \to H^*(S_2)$
in question is therefore the upper route in the diagram.  However, the
map $H^*(\hat{S}_1) \to H^*(S_2)$ induced by the lower route is a
component of $\hat{d}_1$, the differential for $\hat{M}_g$, which
vanishes since $\hat{f}$ is perfect.  Hence the map $\V \to H^*(S_2)$
must vanish.

For the second component, the map $\tau \V \to
H^*(S_2)$, this argument no longer works.  Instead, we resort to a
parity argument.

There is a symplectomorphism $\rho: M_g \to M_g$ induced by a
half-twist of the $g$th handle of the surface $\rs$.
Explicitly, it is given by $A_g \mapsto A_g B_g A_g^{-1} B_g^{-1}
A_g^{-1}$ and $B_g \mapsto A_g B_g^{-1} A_g^{-1}$, with all other
coordinates remaining fixed.  This fixes $f$ and
hence preserves all of the sets discussed above.  It acts as an
involution on $S_1$, preserving the projection to $M_{g-1}$ and the
locus where $B_g = I$, but reversing orientation.  It therefore takes
$\tau$ to $-\tau$, so it acts as $-1$ on $\tau \V$.
On the other hand, it fixes $S_2$, so it acts trivially on
$H^*(S_2)$. 

A more subtle point is that $\rho^*$ commutes with the Thom
isomorphism $H^*(U_2, U_1) \to H^*(S_2)$.  This follows from the Morse
lemma provided that $\rho$ preserves the orientation of the negative
normal bundle.  Of course, $\rho$ will not preserve every Riemannian
metric, so it may not preserve the corresponding negative normal
bundle.  But it does commute with the $\Ui$-action: this can be
verified directly, or deduced from the fact that it preserves $f$ and
the symplectic form.  If a suitable metric is chosen---for example,
one equivalent near $S_2$ to a $\Ui$-invariant flat product metric on the
normal bundle, via the equivariant symplectic tubular neighborhood
theorem---then its negative normal bundle is the $-1$-weight space for
the $\Ui$-action on $TM_g|_{S_2}$.  This subbundle, and its natural
orientation coming from the $\Ui$-action, are therefore preserved by
$\rho$.
 
Any element of $\tau \V$ is therefore acted on as $-1$ in the first
entry, and hence in the last entry, of the sequence of natural maps
$$H^*(S_1)  \lrow  H^*(U_1)  \lrow  H^*(U_2, U_1) \lrow H^*(S_2).$$ 
Its image in $H^*(S_2)$ must therefore be zero.

\br{The map \boldmath $H^*(S_2) \to H^*(S_3)$}  
\label{s2s3}
This map is handled exactly like the previous one.  The fact that
$\hat{f}$ is perfect implies the vanishing of the component $H^*(S_2)
\to \tau \V$.  Indeed, there is a diagram 
{\def\arraycolsep{1pt}
$$
\begin{array}{ccccccccccccc}
&&&&&&&& H^*(U'_3, U'_2) & \lrow & H^*(S_3) && \\
&&&&&&& \nearrow &&&& \searrow & \\ 
H^*(S_2) & \lrow & H^*(U_2, U_1) & \lrow & H^*(U_2) & \lrow &
H^*(U'_2) &&&&&&
H^*(M_{g-1}) \\
&&&&&&& \searrow &&&& \nearrow & \\
&&&&&&&& H^*(\hat{U}'_3, \hat{U}'_2) & \lrow &
H^*(\hat{S}_3) 
&&
\end{array}$$}
where $U'_3 = f^{-1}(-\ratio23, \ratio32)$, $\hat{U}'_3 =
\hat{f}^{-1}(-\ratio23, \ratio32)$, and the final column consists of
push-forwards.  It commutes because the map $H^*(U'_2) \to H^*(S_3)$
can be interpreted as restriction to $f^{-1}(\ratio13)$, which is a
sphere bundle over $S_3$, followed by push-forward by projection, and
similarly for $H^*(U'_2) \to H^*(\hat{S}_3)$.  But the push-forward
$H^*(S_3) \to H^*(M_{g-1})$ is projection on the factor $\tau \V$.

On the other hand, the parity argument implies the vanishing of the
component $H^*(S_2) \to \V$.  Indeed, we have already
seen that $\rho^*$ commutes with all of the arrows between $H^*(S_2)$
and $H^*(U_3,U_2)$ in the above diagram.  However, the Thom
isomorphism $H^*(U_3, U_2) \to H^*(S_3)$ does not commute with
$\rho^*$.  Rather, the two composites differ by a sign.  This is
because $\rho$ reverses the orientation of $S_3$, but not of the open
set $U_3$; after all, it is a symplectomorphism.  It therefore
reverses the orientation of the normal bundle to $S_3$, so it acts as
$-1$ on the Thom class.  Since any class in $H^*(S_2)$ is acted on
trivially by $\rho$, its image in $H^*(S_3)$ must therefore be acted
on as $-1$, so its component in $\V$ must vanish. \er

\br{The map \boldmath $H^*(S_1) \to H^*(S_3)$}
\label{s1s3} 
At this point we have shown that the first differential vanishes.  The
second differential is therefore the map $H^*(S_1) \to H^*(S_3)$ given
by taking $j=1$ in the description of \S\ref{spec}, and it suffices to
show that this also vanishes.  This time, since both $H^*(S_1)$ and
$H^*(S_3)$ split as $\V \oplus \tau \V$, the differential splits into
four components.

Of these, the component $\V \to \tau \V$ vanishes by the perfect
argument.  Simply note that, since $H^*(U_3, U_2) = H^*(U'_3,
U'_2)$ by excision, a prime may be added to every set in the middle
row of (\ref{diag}) without changing the differential.  Then graft on the
hexagonal diagrams of (\ref{s1s2}) and (\ref{s2s3}) and argue as before.

On the other hand, the component $\V \to \V$ vanishes by the parity
argument.  Indeed, suppose that a class in $H^*(S_1)$ belongs to the
component $\V$, so that it is invariant under $\rho^*$.  Then the same
is true of its image in $H^*(U_1)$.  Since the differential $H^*(S_1)
\to H^*(S_2)$ vanishes, this lifts to a class $u \in H^*(U_2)$.  Now
$u$ may not be invariant under $\rho^*$, but the equivariance of the
middle row of (\ref{diag}) implies that $\rho^* u = u + v$, where $v$
is a class in the image of $H^*(U_2, U_1)$.  Since the differential
$H^*(S_2) \to H^*(S_3)$ vanishes, the image of $v$ in $H^*(U_3, U_2)$
must be zero.  Hence $u$ maps to a class in $H^*(U_3, U_2)$ invariant
under $\rho^*$.  Its image in $H^*(S_3)$ is therefore anti-invariant
as in (\ref{s2s3}), so its component in $\V$ vanishes.

Likewise, the component $\tau \V \to \tau \V$ vanishes by the same
parity argument, with the roles of invariance and anti-invariance
interchanged.

This leaves only one of the four components, namely the map $\tau \V
\to \V$.  Now $\tau$ is the restriction to $S_1$ of a global class
$\tau \in H^3(M_g)$, namely the Poincar\'e dual of the locus where
$B_g = I$.  Furthermore, the differential $d_2$ is a module
homomorphism over $H^*(M_g)$.  Indeed, multiplication by a class such
as $\tau$ clearly commutes with the natural maps shown as horizontal
and vertical arrows in (\ref{diag}).  It is therefore compatible with
the lifting from $H^*(U_1)$ to $H^*(U_2)$: given a lifting $u \in
H^*(U_2)$ of a class $w \in H^*(U_1)$, $\tau u$ is a lifting of $\tau
w$.  And it is compatible with the Thom isomorphisms, since they are
simply given by the cup product with the Thom class.  Consequently,
for any $v \in V$, $d_2(\tau v) = \tau d_2(v)$, but $d_2(v)$ has
already been shown to vanish.  This completes the proof of the main
theorem. \er

\bit{The case of additional punctures}
\label{punct}

As promised in the introduction, the main theorem goes through for
flat $\Suii$ connections on a surface $X$ with more than one puncture.
The proof is essentially the same.  However, it is necessary to
generalize the background material of \S\ref{bkg} to the case of
additional punctures.

So let $p_1, \dots, p_n$ be distinct points in $\rs$, let $t_1, \dots,
t_n \in [0,1]$, let $\cl_1, \dots, \cl_n \subset \Suii$ be the
conjugacy classes in $\Suii$ containing 
$\diag(e^{i \pi t_j}, e^{-i \pi t_j})$, 
and let $M_{g,n}$ denote the moduli space of flat
$\Suii$ connections on $\rs$, or rather on $\rs \sans \{ p_1, \dots,
p_n \}$, having holonomy around $p_j$ in $\cl_j$.
Of course $M_{g,n}$ depends on the choice of the $t_j$.

In analogy with \S\ref{bkg}, $M_{g,n}$ can be described as
follows.  Let 
$\mu_{g,n}: \Suii^{2g} \times \prod_{j=1}^n \cl_j \to \Suii$
be given by  
$$\mu_{g,n}(A_1, B_1, A_2, B_2, \dots, A_g, B_g, C_1, \dots, C_n) = 
\Left( \prod_{i=1}^g [A_i, B_i] \, \Right)  
\Left( \prod_{j=1}^n C_j \Right).$$
Then $M_{g,n} = \mu_{g,n}^{-1}(I) / \Suii$.

\bs{Proposition}%
\label{critval}%
For $J \subset \{1, \dots, n\}$, let $\kappa_J = \half \Left( \sum_{j
\in J} t_j - \sum_{j \not\in J} t_j \Right)$.  Then $\mu_{g,n}$
has $I \in \Suii$ as a critical value if and only if for some $J$,
$\kappa_J \in \Z$.   \es

\pf.  The derivative of $\mu_{g,n}$ at $(A_i, B_i, C_j)$ is a map 
$D \mu_{g,n} : \suii^{2g} \oplus
\bigoplus_j T_{C_j} \cl_j \to \suii$.
A direct calculation shows that 
\begin{eqnarray*}
D\mu_{g,n}(a_i, b_i, c_j) 
& = & \ad \Left( \prod_{j=1}^n C_j \Right)^{-1}
D \mu_{g,0}(a_i, b_i) + D \mu_{0,n}(c_j) \\
& = & \sum_{i=1}^g 
\ad \Left( \prod_{k>i} [A_k, B_k] \prod_{j=1}^n C_j \Right)^{-1}
\ad B_i A_i 
\Left( (\ad B_i^{-1} -1) \, a_i + (1-\ad A_i^{-1}) \, b_i \Right) \\
&  & \phantom{xxx} 
+ \sum_{j=1}^n \ad \Left( \prod_{\ell>j} C_\ell \Right)^{-1} \, c_j.
\end{eqnarray*}

Suppose first that $g=0$.  Any $c_j \in T_{C_j}\cl_j$ is of
the form $(1 - \ad C_j^{-1}) \, d_j$ where $d_j \in \suii$.  A further
computation shows that
$$D \mu_{0,n} (c_j) = \sum_{j=1}^{n-1} \left( \ad \Left(\prod_{k>j}
    C_k\Right)^{-1} - 1\right) (d_{j+1} - d_j).$$
This fails to be surjective if and only if for all $j$, the image of
$\ad \left( \prod_{k>j} C_k \right)^{-1}-1$ is the same, which holds
if and only if all $\prod_{k>j} C_k$ commute, which holds if and only
if all $C_j$ commute, since $\prod_{j=1}^n C_j = I$.  The $C_j$ may
then be simultaneously diagonalized to \ $\diag(e^{\pm i \pi t_j},
e^{\mp i \pi t_j})$.  Let $J$ be the set of $j$ for which the plus
sign holds in the first factor.  Then the stated condition on the
$t_j$ holds.

For $g>0$, the presence of the term 
$(\ad B_i^{-1} -1) \, a_i + (1-\ad A_i^{-1}) \, b_i$
in the formula above shows that if $D \mu_{g,n}$ is not surjective,
then $A_i$ and $B_i$ must commute for each $i$.  This brings us back
to the case $g=0$.  \fp

The case studied before was that of one puncture with $t_1 = 1$.
Without loss of generality we may assume that, if we are not in this
case, then each $t_j \in (0,1)$.  Indeed, those $j$ with $t_j = 0$ or
1 can be eliminated from the construction above without changing it,
unless there are an odd number of $j$ with $t_j = 1$.  In the latter
case one can multiply some other $C_j$, say $C_1$, by $-I$, and change
$t_1$ to $1-t_1$.

Suppose that the $t_j$ are so chosen that Proposition \ref{critval} is
false, and $I$ is a regular value.  Then $M_{g,n}$ is smooth.  A
symplectic form on $M_{g,n}$ has been constructed by Guruprasad et
al.\ \cite{guru}.  Moreover, if $f: M_{g,n} \to [-1,1]$ and a
$\Ui$-action on $f^{-1}(-1,1)$ are defined just as in \S\ref{bkg},
ignoring the extra $C_j$ variables, then Audin shows \cite{audin2} that
this action is again symplectic with moment map $-i\arccos f$.

\br{The critical submanifolds of \boldmath$f$}%
\label{crit2}%
The critical submanifolds of $f$ are classified just as in
(\ref{crit}).  First, there are $f^{-1}(-1)$ and $f^{-1}(1)$, both of
which are $\Suii$-bundles over $M_{g-1,n}$ with vanishing Euler class.
These are actually empty in certain cases when $g=1$, but the
main theorem will still hold even then. 

The remaining critical points are again fixed points of the
$\Ui$-action.  It is straightforward to check that these are all
conjugate to $2g+n$-tuples such that
$$A_g = \left(\begin{array}{cc}e^{-i\pi\kappa}&0\\0&e^{i\pi\kappa}\
\end{array}\right), \, \, \,
B_g = \left(\begin{array}{rc}0&1\\-1&0\end{array}\right),$$ 
and the remaining $A_i$, $B_i$, and $C_j$ are all diagonal.  
For $i<g$, the $A_i$ and
$B_i$ are free to move in the 1-parameter subgroup of diagonal
matrices, but each $C_j$ is in a fixed conjugacy class
distinct from $\pm I$ and so must equal 
$\diag(e^{\pm i\pi t_j}, e^{\mp i\pi t_j})$.  Let $J$ be the set of
those $j$ for which the upper sign holds.  
The constraint imposed by $\mu$ then implies that 
$\kappa = \kappa_J$, as defined in Proposition \ref{critval}.
The critical submanifolds in $f^{-1}(-1,1)$ therefore consist of a
disjoint union of $2^n$ tori of dimension $2g-2$.
\er

\bs{Theorem}
The map $f$ is a perfect Bott-Morse function on $M_{g,n}$.
\es

\pf.  The proof is similar to that of the main theorem.  We indicate
only the points where it must be modified.

First of all, the three lemmas at the beginning of \S\ref{pf} go
through unchanged, except that the formula in the proof of Lemma
\ref{lemma} gets replaced by
$$D \mu_{g,n} = D \mu_{g-1,n} + \ad \Left( \prod_j C_j \Right)^{-1} 
D \mu_{1,0}.$$

Suppose that $f^{-1}(\pm 1) \neq \emptyset$, since
the result follows directly from Theorem \ref{kir} otherwise.  Write
the set of critical points as a disjoint union of submanifolds $S_1,
\dots, S_N$ so that
$$-1 = f(S_1) < f(S_2) < \cdots < f(S_{N-1}) < f(S_N) = 1.$$
If there are no coincidences between the values of $\kappa_J$ obtained
in (\ref{crit2}) for different $J$, then $N$ will equal $2^n + 2$.
Otherwise some of the $S_i$ will be disjoint unions of tori.  The
Morse indices of the components may vary, so the Thom class will be
inhomogeneous, but this makes no difference.

Since there are $N$ critical values, only the first $N-1$
differentials in the spectral sequence can possibly be nonzero.  By
induction suppose all differentials before $d_k$ vanish.  If $1 < i <
k+i <n$, then the component $H^*(S_i) \to H^*(S_{k+i})$ of $d_k$
vanishes because $-i \arccos f$ is a moment map there.  The component
$H^*(S_1) \to H^*(S_{k+1})$ is handled as in (\ref{s1s2}), except that
in the parity argument, the $\rho$-equivariance of the long exact
sequence of the relative cohomology of the pair $(U_{i+1}, U_i)$ must
be applied inductively as in (\ref{s1s3}) to show that the
differential preserves the $-1$-weight space.  Likewise the component
$H^*(S_{N-k}) \to H^*(S_N)$ is handled as in (\ref{s2s3}).  Finally,
if $k=N-1$, then $d_k: H^*(S_1) \to H^*(S_N)$ is handled as in
(\ref{s1s3}).  \fp

This leads to a formula for the Betti numbers of $M_{g,n}$ as soon as
the Morse indices of the critical tori are computed.  At this point
our approach meets the route taken by Hoyle \cite{hoy}.  From a
lengthy trigonometric computation, he deduced the following formula.

\bs{Lemma (Hoyle)}
For any $J \subset \{1, \dots, n \}$, the index of the critical torus
associated to $J$ is $2g+2n-2|J|+4[\kappa_J]$.
\es

Therefore, the Poincar\'e polynomials of $M_{g,n}$ satisfy the recursion 
\begin{eqnarray*}
P_t(M_{g,n}) 
& = & P_t(S_1) + t^3 P_t(S^3) + (1+t)^{2g-2} 
\sum_J t^{2(g+n-|J|+2[\kappa_J])} \\
& = & (1+t^3)^2 P_t(M_{g-1,n}) + (t+t^2)^{2g-2}
\sum_J t^{2(n+1-|J|+2[\kappa_J])}.
\end{eqnarray*}
This immediately implies the following.

\bs{Corollary}
The Poincar\'e polynomial of $M_{g,n}$ satisfies
$$P_t(M_{g,n}) = (1+t^3)^{2g} P_t(M_{0,n}) + \frac{(1+t^3)^{2g} -
  (t+t^2)^{2g}}{(1-t^2)(1-t^4)} \sum_J t^{2(n+1-|J|+2[\kappa_J])}.$$
\es
\vspace{-3ex}
From this, one can obtain an explicit expression for $P_t(M_{g,n})$ by
substituting Hoyle's formula for $P_t(M_{0,n})$.  
  
\br{Remark} Assuming the Harder-Narasimhan formula, the recursion in
the genus could also be proved by algebro-geometric methods.  Here is
a sketch of the argument.  By a theorem of Mehta-Seshadri \cite{ms},
$M_{g,n}$ may be regarded as a moduli space of vector bundles over $X$
with parabolic structure at each $p_j$.  If $t_n = 1$ but the
remaining $t_j$ are small, then $M_{g,n}$ is a
$(\C\Pj^1)^{n-1}$-bundle over $M_g$ for $g >0$, so the recursion holds
in this case.  For general $t_j$, choose a path in $[0,1]^n$
connecting $(t_1, \dots t_n)$ to the previous special case in such a
way that for any point on the path, $\kappa_J \in \Z$ for at most one
$J \subset \{ 1, \dots, n \}$.  Then the moduli spaces on either side
of each such point are related by a blow-up and blow-down centered on
a disjoint union of tori \cite{bh, flip}.  It can be checked that this
alters the Betti numbers in a way which preserves the recursion.
Unlike the Morse approach, however, this algebraic approach gives no
insight about why the recursion holds.  \er

\bit{Some final remarks}
\label{action}

Many other cheerful facts about $M_{g,n}$ can be deduced from the
argument used to prove the main theorem.  Let us mention two of them.

\bs{Proposition}
The integral cohomology of $M_{g,n}$ is torsion-free.
\es

\pf.  When $g=0$, then Hoyle \cite{hoy} exhibits a Hamiltonian circle
action on $M_{0,n}$ whose fixed points are isolated except for
two copies of $M_{0,n-1}$.  By induction we may suppose these are
torsion-free, and Frankel \cite[Cor.\ 1]{frank} then shows that
$M_{0,n}$ is torsion-free.

Now suppose $g>0$ and consider again the Bott-Morse function $f$.  The
spectral sequence of $\S2$ works equally well with integral
cohomology.  Everything in the proof of the main theorem goes through,
except that the differentials over $\Z$ might take nonzero values in
the torsion part of $H^*(S_2, \Z)$ or $H^*(S_3, \Z)$.  But $S_2$ is a
torus, so it is certainly torsion-free, and $S_1$ and $S_3$ may be
assumed torsion-free by induction on $g$.  \fp

\bs{Proposition}
Let $(\Z / 2)^{2g}$ act on $M_{g,n}$ by 
$$(\delta_i, \epsilon_i) \cdot
(A_i, B_i, C_j) = ((-1)^{\delta_i} A_i, (-1)^{\epsilon_i} B_i, C_j).$$
Then the induced action on $H^*(M_{g,n}, \Q)$ is trivial.  \es

\pf.  The action of $(\Z/2)^{2g-2}$ on the first $2g-2$ factors
preserves the $\Ui$-action, the map $f$ and so on.  The whole proof of
the main theorem can therefore be $(\Z/2)^{2g-2}$-graded. But since
$(\Z/2)^{2g-2}$ acts trivially on $\tau$ and on $H^*(S_2)$, by
induction the whole grading is trivial.  Hence $(\Z/2)^{2g-2}$ acts
trivially on $H^*(M_{g,n})$.

The last two factors are really no harder.  After all, the choice of
an ordering on the handles of $\rs$ was arbitrary.  For example, the
argument still works if $f$ is replaced by $\half \tr A_1$, and so on.
\fp

The value of this last result is to relate the cohomology of
$M_{g,n}$ to that of the corresponding moduli space of flat $\Uii$
connections.  Let $\tilde{M}_{g,n}$ be the moduli space of flat $\Uii$
connections on $X \sans \{ p_1, \dots, p_n \}$ with holonomy around
$p_j$ in the conjugacy class $\cl_j$.  Just as $M_{g,n}$ was,
$\tilde{M}_{g,n}$ may be described in terms of a map
$\tilde{\mu}_{g,n}$.

\bs{Corollary}
As rings, $H^*(\tilde{M}_{g,n}, \Q) \cong H^*(\Ui^{2g}, \Q) \otimes
H^*(M_{g,n}, \Q)$.   
\es

\pf.  There is a
natural map $\Ui^{2g} \times M_{g,n} \to \tilde{M}_{g,n}$ given by 
$$(\kappa_i, \lambda_i), (A_i, B_i, C_j) \mapsto (\kappa_i A_i, \la_i
B_i, C_j).$$
Indeed, $\tilde{M}_{g,n} = (\Ui^{2g} \times M_{g,n}) /
(\Z/2)^{2g}$, where $(\Z/2)^{2g}$ acts diagonally on $\Ui^{2g}$ and
$M_{g,n}$ as above.  The induced action on $H^*(\Ui^{2g}, \Q)$ is
certainly trivial, and the induced action on $H^*(M_{g,n}, \Q)$ is
trivial by the lemma above.  The proof is completed by the result of
Grothendieck \cite{toh} that the rational cohomology ring of a
quotient by a finite group is the invariant part of the rational
cohomology.  \fp

Compare Newstead \cite{n}, Harder-Narasimhan \cite{hn} 
and Atiyah-Bott \cite{ab} in the case $n=0$.

\end{document}